\documentstyle{article}

\def\qmod#1#2{{\raise1pt\hbox{$#1$}\kern-1pt\big/
               \kern-1pt\raise-1pt\hbox{$#2$}}}

\newfam\msbfam
\font\tenmsb=msbm10 at 10 pt
\font\sevenmsb=msbm10 at 7pt
\font\fivemsb=msbm10 at 5pt

\textfont\msbfam=\tenmsb
\scriptfont\msbfam=\sevenmsb
\scriptscriptfont\msbfam=\fivemsb

\def\Bbb{\fam\msbfam\tenmsb}

\def\C{{\Bbb C}}

\def\H{{\Bbb H}}

\def\P{{\Bbb P}}

\def\R{{\Bbb R}}
\def\Z{{\Bbb Z}}

\newfam\cmbfam
\font\tencmb=cmbx10
\font\sevencmb=cmbx7
\font\fivecmb=cmbx5

\textfont\cmbfam=\tencmb
\scriptfont\cmbfam=\sevencmb
\scriptscriptfont\cmbfam=\fivecmb

\def\map{\longrightarrow}
\def\textmap#1{\mathop{\vbox{\ialign{
                                ##\crcr
    ${\scriptstyle\hfil\;\;#1\;\;\hfil}$\crcr
    \noalign{\kern-1pt\nointerlineskip}
    \rightarrowfill\crcr}}\;}}

\def\textlmap#1{\mathop{\vbox{\ialign{
                                ##\crcr
    ${\scriptstyle\hfil\;\;#1\;\;\hfil}$\crcr
    \noalign{\kern-1pt\nointerlineskip}
    \leftarrowfill\crcr}}\;}}

\newfam\meuffam
\font\tenmeuf=eufm10
\font\sevenmeuf=eufm7
\font\fivemeuf=eufm5

\textfont\meuffam=\tenmeuf
\scriptfont\meuffam=\sevenmeuf
\scriptscriptfont\meuffam=\fivemeuf

\def\germ{\fam\meuffam\tenmeuf}

\def\cg{{\germ c}}

\def\picture#1by#2(#3){
\vbox to #2 {
  \hrule width #1 height 0pt depth 0pt \vfill \special{picture #3}}
}

\def\scaledpicture#1by#2(#3scaled#4){{
\dimen0=#1  \dimen1=#2
\divide\dimen0 by 1000 \multiply\dimen0 by #4
\divide\dimen1 by 1000 \multiply\dimen1 by #4
\picture \dimen0 by \dimen1 (#3 scaled #4)}}
\def\dfigure#1by#2(#3scaled#4offset#5:#6)
  {\medskip
   \vglue 2mm minus 2mm
   $$
     \hbox{
       \hglue#5
       {\scaledpicture #1 by #2 (#3 scaled #4)}
     }
   $$
   \par\goodbreak
   \vglue 2mm minus 2mm
   \medskip}

\begin{document}
\def\Pr{{\rm Pr}}
\def\tr{{\rm Tr}}
\def\ad{{\rm ad}}
\def\End{{\rm End}}
\def\Pic{{\rm Pic}}
\def\NS{{\rm NS}}
\def\deg{{\rm deg}}
\def\Hom{{\rm Hom}}
\def\Aut{{\rm Aut}}
\def\Herm{{\rm Herm}}
\def\Vol{{\rm Vol}}
\def\pf{{\bf Proof: }}
\def\id{{\rm id}}
\def\im{{\rm im}}
\def\rk{{\rm rk}}
\def\coker{{\rm coker}}
\def\Sp{{\rm Sp}}
\def\Spin{{\rm Spin}}
\def\h{{\bf H}}
\def\dv{\bar\partial}
\def\dva{\bar\partial_A}
\def\da{\partial_A}
\def\p{\partial\bar\partial}
\def\pa{\partial_A\bar\partial_A}
\def\Dr{{\raisebox{.16ex}{$\not$}}\hskip -0.35mm{D}}
\def\oo{{\scriptstyle{\cal O}}}
\newtheorem{sz}{Satz}
\newtheorem{szfr}{Satzfr}
\newtheorem{thry}[sz]{Theorem}
\newtheorem{thfr}[szfr]{Th\'eor\`eme}
\newtheorem{pr}[sz]{Proposition}
\newtheorem{re}[sz]{Remark}
\newtheorem{co}[sz]{Corollary}
\newtheorem{cofr}[szfr]{Corollaire}
\newtheorem{dt}[sz]{Definition}
\newtheorem{lm}[sz]{Lemma}

\noindent { }
\\ \\
\centerline{\Large{\bf Seiberg-Witten invariants for manifolds
with $b_+=1$}}
\\  \\
\vspace*{1.5mm}
\centerline{\large Christian Okonek\footnotemark  \ \ and \ \
Andrei
Teleman$^1$}\\ \footnotetext[1]{Partially supported by:
AGE-Algebraic Geometry in
Europe, contract No ERBCHRXCT940557 (BBW 93.0187), and by  SNF,
 nr. 21-36111.92}
{\footnotesize{\sl Abstract --}} {\footnotesize In this note we
describe the
Seiberg-Witten invariants, which have been introduced in [W],
for manifolds with
$b_+=1$.  In this case the invariants depend on a chamber
structure, and there
exists a universal wall crossing formula. For every K\"ahler
surface with $p_g=0$
and $q$=0, these invariants are non-trivial for all
$\Spin^c(4)$-structures of
non-negative index. }
\vskip 4pt

\centerline{{\bf Les invariants de Seiberg-Witten pour les
vari\'et\'es avec $b_+=1$ }}
\vskip 4pt
{\footnotesize{\sl R\'esum\'e --}} {\footnotesize   Dans cette note
nous d\'ecrivons les invariants de Seiberg-Witten introduits dans
[W], pour les vari\'et\'es telles que $b_+=1$. Ces invariants
d\'ependent d'un param\`etre auxiliaire variant dans l'ensemble
des chambres, et il existe une formule universelle de passage
\`a travers un mur. Pour les surfaces k\"ahl\'eriennes telles que
$p_g=0$ et $q=0$, les invariants associ\'es \`a toute
$\Spin^c(4)$-structure d'index non-n\'egatif sont non-triviaux.  }
\vskip 8pt

{\sl Version fran\c{c}aise abr\'eg\'ee } - Soit $X$ une 4-vari\'et\'e
ferm\'ee orient\'ee connexe, et $c\in H^2(X,\Z)$ tel que
$c\equiv w_2(X)$ (mod 2).  Soit $\hat P$ un $\Spin^c(4)$-fibr\'e
\underbar{compatible}, c'est-\`a-dire que $c_1(\hat P\times_{\det}\C)=c$
et la $GL_+(4,\R)$-extension $\hat P\times_{\tilde\pi}GL_+(4,\R)$  est
isomorphe au fibr\'e des rep\`eres orient\'es de $\Lambda^1_X$.
Soient $\Sigma^{\pm}:=\hat P\times_{\sigma_{\pm}}\C^2$
les fibr\'es des spineurs associ\'es, de sorte que $\det\Sigma^{\pm}=\hat
P\times_{\det}\C$ [OT1].
Une \underbar{application}  \underbar {de} \underbar{Clifford}
du type $\hat P$ est un $GL_+(4,\R)$-isomorphisme
 $\gamma:\Lambda^1_X\rightarrow
P\times_{\pi}\R^4$.  Une telle application $\gamma$ d\'efinit
une m\'etrique $g_\gamma $,  et des isomorphismes
$\Gamma:\Lambda^2_{\pm}\map su(\Sigma^{\pm})$.   D\'esignant
 par
${\cal C}={\cal C}(\hat P)$ l'espace  des applications de Clifford
du type $\hat P$, on a un isomorphisme naturel
$\qmod{{\cal C}}{ \Aut( \hat P)}
\textmap{\simeq}{\cal M}et_X\times
\pi_0\left(\qmod{{\cal C}} { \Aut( \hat P)} \right)$,
o\`u le deuxi\`eme  facteur param\'etrise l'ensemble des classes
d'\'equivalence des $\Spin^c(4)$-structures de classe de Chern  $c$
sur $(X,g)$, pour toute  m\'etrique $g$.
On d\'esigne par $\cg_\gamma\in\pi_0(\qmod{{\cal C}}{ \Aut( \hat
P)})$  l'\'el\'ement d\'efini par $[\gamma]\in \qmod{{\cal C}}{
\Aut( \hat P)})$.
La donn\'ee d'une application de Clifford $\gamma$ d\'efinit une
bijection entre les connexions unitaires sur
$\hat P\times_{\det}\C$ et les $\Spin^c(4)$-connexions sur $\hat P$
qui se projettent (via $\gamma$) sur la connexion de Levi-Civita
de ${g_\gamma}$; elle permet
d'associer \`a toute connexion $A\in{\cal A}(\hat P\times_{\det}\C)$
un op\'erateur de  Dirac  $\Dr_A$.
\vskip 5pt
{\sc D\'efinition.} - {\it Soit $\gamma$ une application de Clifford,
et soit
$\beta\in Z^2_{\rm DR}(X)$ une 2-forme ferm\'ee. Les \'equations
$\beta$-translat\'ees de Seiberg-Witten  s'\'ecrivent
$$ \Dr_{A}\Psi = 0\  ,\ \
\Gamma\left((F_A+{2\pi i}\beta)^+\right) = 2(\Psi\bar\Psi)_0\ .
\eqno{(SW^{\gamma}_\beta)}$$
}
Soit ${\cal W}_{X,\beta}^{\gamma}$ l'espace des modules des
solutions $(A,\Psi)\in{\cal A}(\det \Sigma^+)\times A^0(\Sigma^+)$
de $(SW^{\gamma}_\beta)$ modulo l'action naturelle
du groupe de jauge ${\cal G}={\cal C}^{\infty}(X,S^1)$.
L'espace ${\cal W}_{X,\beta}^{\gamma}$ ne d\'epend  que de
$(g_\gamma,\cg_\gamma)$ et $\beta$, \`a isomorphisme
\underbar{canonique} pr\`es.

L'espace
${\cal B}(c)^*:= \qmod{{\cal A}(\det\Sigma^+)
\times(A^0(\Sigma^+)\setminus\{0\})}{{\cal G}}$ a le type
d'homotopie faible de $K(\Z,2)\times
K(H^1(X,\Z),1)$; il y a un isomorphisme naturel $\nu:\Z[u]\otimes
\Lambda^*(\qmod{H^1(X,\Z)}{\rm Tors})\rightarrow
 H^*({\cal B}(c)^*,\Z)$.

 Soit $c\in H^2(X,\Z)$ un \'el\'ement caract\'eristique. Une
paire
$(g,b)\in{\cal M}et_X\times H^2_{\rm DR}(X)$ est dite
 $c$-\underbar{bonne} si le repr\'esentant
$g$-harmonique de $(c-b)$ n'est pas anti-autodual.

Soit maintenant   $(g,b)$ une paire c-bonne  et
$\cg\in\pi_0(\qmod{{\cal
C}}{\Aut(\hat P)})$ . L'espace
${\cal W}^{\gamma}_{X,\beta}$ est une vari\'et\'e lisse compacte
de dimension
$w_c:=\frac{1}{4}(c^2-2e(X)-3(\sigma(X))$  pour toute
$\gamma\in{\cal C}$ qui se projette sur
$(g,\cg)$ et tout $\beta$ dans un sous-ensemble ouvert
dense de $b$. Une telle forme
$\beta$ sera appel\'ee r\'egul\`ere. La vari\'et\'e
${\cal W}^{\gamma}_{X,\beta}$
peut alors \^etre orient\'ee  en choisissant une orientation
$\oo$ du fibr\'e
$\det(H^1(X,\R))\otimes\det(\H^2_{+,g}(X)^{\vee})$. Soit
$[{\cal
W}^{\gamma}_{X,\beta}  ]_{{\raisebox{-.5ex}{$\oo$}}}
\in H_{w_c}({\cal B}(c)^*,\Z)$
la classe fondamentale associ\'ee  \`a $\oo$.

 La \underbar{forme}  {de}
\underbar{Seiberg}-\underbar{Witten}  associ\'ee \`a la donn\'ee
$({\oo},(g,b),\cg)$ est l'\'el\'ement
 $SW_{X,\oo}^{(g,b)}(\cg)\in \Lambda^*
H^1(X,\Z)$ d\'efini  par
$SW_{X,\oo}^{(g,b)}(\cg)(l_1\wedge\dots\wedge l_r)=
 \left\langle\nu(l_1)\cup\dots\cup\nu(l_r)\cup u^{\frac{w_c-r}{2}},
[{\cal W}^{\gamma}_{X,\beta}   ]_{{\raisebox{-.5ex}{$\oo$}}}
\right\rangle\
$
sur les \'el\'ements d\'ecomposables $l_1\wedge\dots\wedge l_r$
avec $r\equiv
w_c$ (mod 2). Ici $\gamma\in{\cal C}$ se projette sur la paire
 $(g,\cg)$ et $\beta\in b$ est r\'eguli\`ere.
$SW_{X,\oo}^{(g,b)}(\cg)$ ne d\'epend pas du choix de
$\gamma$ et $\beta$.

Si $b_+>1$, $SW_{X,\oo}^{(g,b)}(\cg)$ ne d\'epend pas m\^eme
du choix de la paire $c$-bonne
$(g,b)$, donc on peut d\'esigner cet invariant simplement par
$SW_{X,\oo}(\cg)\in\Lambda^* H^1(X,\Z)$.
\def\ooo{{\scriptscriptstyle{\cal O}}}
Si $b_1=0$, on obtient des nombres que
nous  d\'esignons par  $n_{\cg}^{\ooo}$; ces nombres peuvent \^etre
consid\'er\'es comme des raffinements des
nombres $n_c^{\ooo}$ introduit en [W].  En effet
$n_c^{\ooo}=\sum_{\cg}n_{\cg}^{\ooo}$, o\`u la somme est faite
par rapport \`a
$\cg\in\pi_0(\qmod{{\cal C}}{\Aut(\hat P)})$.

Supposont maintenant que $b_+=1$.  Il y a une application naturelle
${\cal M}et_X\map
\P(H^2_{\rm DR}(X))$ qui associe \`a  $g\in{\cal M}et_X$ la droite
$\R[\omega_+]\subset H^2_{\rm DR}(X)$, o\`u $\omega_+$ est une
2-forme $g$-autoduale harmonique non-triviale. L'espace hyperbolique
${\bf H}:=\{\omega\in H^2_{\rm DR}(X)|\ \omega^2=1\}
$
a deux composantes connexes, et le choix d'une composante
${\bf H}_0$ d\'efinit
une orientation  de la droite   $\H^2_{+,g}(X)$ pour toute
m\'etrique  $g$. Soit $\omega_g$ le g\'en\'erateur de $\H^2_{+,g}(X)$
tel que $[\omega_g]\in {\bf H}_0$ .
\vskip 5pt
{\sc D\'efinition.} - {\it Soit $X$ une 4-vari\'et\'e avec $b_+=1$, et
$c\in H^2(X,\Z)$ un
\'el\'ement caract\'eri\-stique. Le \underbar{mur} associ\'e  \`a $c$
est l'hypersurface
$c^{\bot}:=\{(\omega,b)\in{\bf H}\times H^2_{\rm DR}(X)|\ (c-b)
\cdot\omega=0\}$.
Les composantes connexes de ${\bf H}\setminus c^{\bot}$
seront appel\'ees \underbar{chambres} du type $c$.
}
\vskip 3pt

Les murs ne sont pas  lin\'eaires! Tout \'el\'ement caract\'eristique
$c$ d\'efinit pr\'ecis\'ement quatre chambres du type $c$ :
$C_{{\bf H}_0,\pm}:=\{(\omega,b)\in{\bf H}_0\times
 H^2_{\rm DR}(X)|\
\pm(c-b)\cdot\omega<0\}\ , $
o\`u ${\bf H}_0$ est l'une des deux composantes connexes de
${\bf H}$. Toute chambre contient des paires de la
forme  $([\omega_g],b)$. Choisissons maintenant une orientation
$\oo_1$ de $H^1(X,\R)$.
\vskip 5pt
{\sc D\'efinition.} -  {\it L' invariant de Seiberg-Witten  associ\'e \`a
la donn\'ee $(\oo_1,{\bf H}_0,\cg)$ est la fonction
$
SW_{X,(\oo_1,{\bf H}_0)}(\cg):\{\pm\}\rightarrow
\Lambda^* H^1(X,\Z)
$
d\'efinie par  $SW_{X,(\oo_1,{\bf H}_0)}(\cg)(\pm):=
SW^{(g,b)}_{X,\oo}(\cg)$, o\`u
$\oo$ est l'orientation induite  par $(\oo_1,{\bf H}_0)$, et
$([\omega_g],b)$ est une paire  appartenant \`a $C_{{\bf H}_0,\pm}$. }
\vskip 3pt

On v\'erifie  facilement les relations
$
SW_{X,(-\oo_1,{\bf H}_0)}(\cg)(\pm)=-SW_{X,(\oo_1, {\bf H}_0)}(\cg)(\pm)$
\ et\hfil\break
$SW_{X,(\oo_1,-{\bf H}_0)}(\cg)(\pm)=
-SW_{X,(\oo_1,{\bf H}_0)}(\cg)(\mp).$
Soit $u_c\in\Lambda^2\left(\qmod{H_1(X,\Z)}{\rm
Tors}\right)$ d\'efini par la formule
$u_c(a,b):=\frac{1}{2}\langle c\cup a\cup b,[X]\rangle $,  pour   $a,b\in
H^1(X,\Z)$.
\vskip 5pt
{\sc Th\'eor\`eme.} -  {\it Supposons  $b^+(X)=1$,  soit $l_{\ooo_1}$
le g\'en\'erateur de $\Lambda^{b_1}(H^1(X,\Z))$ d\'efini par
l'orientation $\oo_1$, et soit $r\geq 0$, $r\equiv w_c$
(mod 2) .  Pour tout $\lambda\in
\Lambda^r\left(\qmod{H_1(X,\Z)}{\rm Tors}\right)$   on a \\
\centerline{$SW_{X,(\oo_1,{\bf H}_0)}(\cg)(+)(\lambda)- SW_{X,(\oo_1,{\bf
H}_0)}(\cg)(-)(\lambda)=\frac{ (-1)^{\left[\frac{b_1-r}{2}\right]} }
{ \left[\frac{b_1-r}{2}\right]!} \langle  \lambda\wedge u_c
^{\left[\frac{b_1-r}{2}\right]},l_{\ooo_1} \rangle
$}
si  $r\leq \min(b_1,w_c)$, et la diff\'erence est nulle dans les autres cas.
}
\vskip 3pt

Soit $(X,g)$ une surface k\"ahl\'erienne munie de sa
$\Spin^c(4)$-structure canonique  et soit $\omega_g$ sa forme de
K\"ahler.  Il y a une bijection naturelle entre les classes de
$\Spin^c(4)$-structures
$\cg$ de classe de Chern  $c$ et les fibr\'es   en droites
$M$ dont la classe de Chern v\'erifie $2c_1(M)-c_1(K_X)=c$.  On
d\'esigne par
$\cg_M$ la classe  d\'efinie par $M$.   Soit  ${\cal D}ou(m)$ l'espace
de Douady  des diviseurs effectifs $D$ sur $X$ tels que
$c_1({\cal O}_X(D))=m$.
\vskip 5pt
{\sc Th\'eor\`eme.} - {\it Soit  $(X,g)$ une surface k\"ahl\'erienne
connexe, et soit $\cg_M$ la classe des
  $\Spin^c(4)$-structures associ\'ee au fibr\'e en droites   $M$ de
classe de Chern
$c_1(M)=m$. Soit $\beta\in A^{1,1}_{\R}$ une forme   repr\'esentant
la classe $b$ telle que
$ (2m-c_1(K_X)-b)\cup[\omega_g]<0$ ($>0$).
\hfill{\break}
i)  Si $c\ \not\in\  NS(X)$, on a ${\cal W}_{X,\beta}^{\gamma_M}
=\emptyset$. Si $c\in NS(X)$, il  existe un isomorphisme r\'eel
analytique naturel
${\cal W}_{X,\beta}^{\gamma_M}\simeq  {\cal D}ou(m)$
$({\cal D}ou(c_1(K_X)-m))$.
 \hfill{\break}
ii)  ${\cal W}_{X,\beta}^{\gamma_M} $ est lisse en un point
correspondant \`a $D\in{\cal D}ou(m)$ si et seulement si
$h^0({\cal O}_D(D))=\dim_D{\cal D}ou(m)\ .$
Cette condition est toujours satisfaite  si $h^1({\cal O}_X)=0$.
\hfill{\break}
iii) Si  ${\cal W}_{X,\beta}^{\gamma_M}$ est lisse en  un point
correspondant \`a $D\in{\cal D}ou(m)$, il a la dimension
pr\'esum\'ee en ce point si et seulement si $h^1({\cal O}_D(D))=0$.
}
\vskip 3pt
Toute surface complexe connexe avec $p_g>0$ est diff\'eomorphe
\`a une surface poss\'edant un diviseur canonique 0-connexe. De l\`a
on d\'eduit une d\'emonstration simple du
\vskip 5pt
{\sc Corollaire.} - {\it ([W] )Tous les invariants de Seiberg-Witten
non-triviaux d'une surface k\"ahl\'erienne avec $p_g>0$ sont d'index 0.}
\vskip 3pt
Au contraire, si $p_g=0$, on a
\vskip 5pt
{\sc Corollaire.} - {\it Soit $X$ une surface k\"ahl\'erienne avec
$p_g=0$ et $q=0$. Pour toute donn\'ee $({\bf H}_0,\cg)$ avec
$w_c\geq 0$, on a
$SW_{X,{\bf H}_0}(\cg)(\{\pm\})=\{0,1\}$ ou $SW_{X,{\bf
H}_0}(\cg)(\{\pm\})=\{0,-1\}$}.
\vskip 3pt
\noindent\hbox to 5cm{\hrulefill} \\

1.  {\sc The twisted  Seiberg-Witten equations}. - Let $X$ be a
closed connected oriented    4-manifold, and let $c\in H^2(X,\Z)$ be a
class with $c\equiv w_2(X)$ (mod 2). A
\underbar{compatible} $\Spin^c(4)$-bundle is a    $\Spin^c(4)$-bundle
$\hat P$
 with $c_1(\hat P\times_{\det}\C)=c$ such that its
$GL_+(4,\R)$-extension
$\hat P\times_{\tilde\pi}GL_+(4,\R)$ is isomorphic to the bundle
of oriented frames in
$\Lambda^1_X$.  Let $\Sigma^{\pm}:= \hat
P\times_{\sigma_{\pm}}\C^2$ be the associated spinor bundles with
$\det\Sigma^{\pm}= \hat P\times_{\det}\C$  [OT1].
\vskip 5pt
{\sc Definition.} - {\it A \underbar{Clifford}   \underbar{map}   of
type $\hat P$ is a $GL_+(4,\R)$-isomorphism  \hbox{
$\gamma:\Lambda^1_X\rightarrow\hat P\times_{\pi}\R^4$}.  }
\\
The $SO(4)$-vector bundle $\hat P\times_\pi\R^4$ can be
identified with the bundle
$\R SU(\Sigma^+,\Sigma^-)$ of real multiples of $\C$-linear
isometries  of determinant 1  from $\Sigma^+$ to $\Sigma^-$.
A Clifford map
$\gamma$ defines a metric
$g_\gamma $ on $X$, a lift $\hat P\map P_{g_\gamma}$ of the
associated frame bundle, and it induces isomorphisms
$\Gamma:\Lambda^2_{\pm}\map su(\Sigma^{\pm})$  [OT1].   Let
${\cal C}={\cal C}(\hat P)$ be the space of all Clifford maps of type
$\hat P$. There is  a  natural  isomorphism
$\qmod{{\cal C}}{ \Aut( \hat P)}
\textmap{\simeq}{\cal M}et_X\times
\pi_0\left(\qmod{{\cal C}} { \Aut( \hat P)} \right)
$, where the second factor is a ${\rm Tors}_2H^2(X,\Z)$-torsor; it
parametrizes the set of equivalence classes of
$\Spin^c(4)$-structures  with Chern class $c$ on $(X,g)$,  for an
arbitrary  metric $g$.
We use the symbol $\cg$ to denote elements in   $\pi_0(\qmod{{\cal
C}}{ \Aut(
\hat P)})$, and we denote by $\cg_\gamma$ the connected
component defined by
$[\gamma]\in\qmod{{\cal C}} { \Aut( \hat P)}$.
A fixed Clifford map $\gamma$ defines a bijection
between unitary connections in
$\hat P\times_{\det}\C$ and $\Spin^c(4)$-connections in $\hat P$
which lift (via
$\gamma$) the Levi-Civita connection in
$P_{g_\gamma}$, and allows to associate  a Dirac operator  $\Dr_A$
to  a connection
$A\in{\cal A}(\hat P\times_{\det}\C)$.
\vskip  5pt
{\sc Definition.} - {\it Let $\gamma$ be a Clifford map, and let
$\beta\in Z^2_{\rm DR}(X)$ be a closed 2-form. The
$\beta$-twisted Seiberg-Witten equations are
$$ \Dr_{A}\Psi = 0\ ,\ \
\Gamma\left((F_A+{2\pi i}\beta)^+\right) = 2(\Psi\bar\Psi)_0\ .
 \eqno{(SW^{\gamma}_\beta)}
$$
}
These twisted Seiberg-Witten equations arise
 naturally in
connection with certain non-abelian monopoles [OT2].  They should
\underbar{not} be regarded as perturbation of $(SW^\gamma_0)$,
since later the cohomology class of $\beta$ will be fixed.

Let ${\cal W}_{X,\beta}^{\gamma} $  be the moduli space  of
solutions $(A,\Psi)\in{\cal A}(\det \Sigma^+)\times A^0(\Sigma^+)$
of
$(SW^{\gamma}_\beta)$ modulo  the natural action
$((A,\Psi),f)\longmapsto
(A^{f^2},f^{-1}\Psi)$ of the gauge group
${\cal G}={\cal C}^{\infty}(X,S^1)$.

The moduli space  ${\cal W}_{X,\beta}^{\gamma}$ depends up to
\underbar{canonical} isomorphism only on
$(g_\gamma,\cg_\gamma)$ and
$\beta$, since  two Clifford maps lifting the same pair $(g,\cg)$
are equivalent modulo $\Aut(\hat P)$.

Now fix a class $b\in H^2_{\rm DR}(X)$, consider
$(SW^{\gamma}_\beta)$ as equation for a triple $(A,\Psi,\beta)\in
{\cal A}(\det \Sigma^+)\times A^0(\Sigma^+)\times b$, and let
${\cal W}_{X,b}^{\gamma}\subset\qmod{  {\cal A}(\det
\Sigma^+)\times A^0(\Sigma^+) \times b }{{\cal G}}$
be the (infinite dimensional) moduli space of solutions.   Finally we
need the universal moduli space
${\cal W}_X \subset\qmod{ {\cal A}(\det \Sigma^+)\times
A^0\ (\Sigma^+)\times Z^2_{DR}(X)\times {\cal C} }
{{\cal G}}$
of solutions of $(SW^{\gamma}_\beta)$ regarded as equations for
tuples
$(A,\Psi,\beta,\gamma)\in{\cal A}(\det \Sigma^+)\times
A^0(\Sigma^+)\times Z^2_{DR}(X)\times {\cal C} $.
We complete the spaces  ${\cal A}(\det\Sigma^+)$,
$A^0(\Sigma^{\pm})$ and $A^2$ with respect to the Sobolev norms
$L^2_q$, $L^2_q$ and $L^2_{q-1}$, and   the gauge group ${\cal G}$
with respect to $L^2_{q+1}$, but we suppress the Sobolev subscripts
in our notations. As usual we denote by the superscript ${\ }^*$ the
open subspace of a moduli space where the spinor component is
non-zero.
\vskip 5pt
{\sc Definition.} - {\it Let $c\in H^2(X,\Z)$ be   characteristic. A pair
$(g,b)\in{\cal M}et_X\times H^2_{\rm DR}(X)$
 is  $c$-good if  the g-harmonic
representant  of $(c-b)$ is not $g$-anti-selfdual.
}
\vskip 5pt
{\sc Proposition.} - {\it Let $X$ be a closed  oriented 4-manifold,
and let $c\in H^2(X,\Z)$ be characteristic. Choose a compatible
$\Spin^c(4)$-bundle $\hat P$  and  an element $\cg\in
\pi_0\left(\qmod{ {\cal C}}{\Aut(
\hat P)} \right)$. \\
i)  The projections
$p:{\cal W}_X \map Z^2_{DR}(X)\times{\cal C}$ and
$p_{\gamma,b}:{\cal W}_{X,b}^{\gamma } \map b$  are proper for all
$\gamma$, $b$.\\
ii)   ${\cal W}_X^* $ and ${{\cal W}_{X,b}^{\gamma}}^* $   are  smooth
manifolds  for all $\gamma$ and $b$.\\
iii)    ${{\cal W}_{X,b}^{\gamma}}^* ={\cal W}_{X,b}^{\gamma} $ if
$(g_\gamma,b)$ is
$c$-good.\\
iv)  If $(g_\gamma,b)$ is $c$-good, then every pair
$(\beta_0,\beta_1)$ of regular values of $p_{\gamma,b}$ can be
joined by a smooth path
$\beta:[0,1]\map b$   such that the fiber product
$[0,1]\times_{(\beta, p_{\gamma,b})}{\cal W}_{X,b}^{\gamma} $
defines a smooth  bordism between   ${\cal W}_{X,\beta_0}^{\gamma}
$ and ${\cal W}_{X,\beta_1}^{\gamma} $. \\
v) If   $(g_0,b_0)$,  $(g_1,b_1)$ are $c$-good pairs which can be
joined by a  smooth path  of $c$-good  pairs,   then there is a
smooth path $(\beta,\gamma):[0,1]\map   Z^2_{\rm DR}(X)
\times{\cal C}$ with the following properties:\\
\hspace*{0.8cm}1. $[\beta_i]=b_i$ and $g_{\gamma_i}=g_i$ for
$i=0,  1$.  \\
\hspace*{0.8cm}2. $\gamma_t$ lifts $(g_{\gamma_t},\cg)$ and
$(g_{\gamma_t},[\beta_t])$ is $c$-good for every $t\in[0,1]$.\\
\hspace*{0.8cm}3.   $[0,1]\times_{((\beta,\gamma),p)}{\cal
W}_X ^*$ is a smooth bordism between
${\cal W}_{X,\beta_0}^{\gamma_0} $ and
${\cal W}_{X,\beta_1}^{\gamma_1} $.\\
vi) If $b_+>1$, then any two $c$-good pairs $(g_0,b_0)$,
$(g_1,b_1)$  can be joined by a smooth path of $c$-good pairs .
}
\vskip 3pt
The proof uses techniques from [DK], [KM]  and [T].
\vskip 3pt

2. {\sc Seiberg-Witten invariants for 4-manifolds with $b_+=1$.} -
Let $X$ be a closed connected   oriented 4-manifold , $c$ a
characteristic element,  and
$\hat P$ a compatible $\Spin^c(4)$-bundle. We put
${\cal B}(c)^*:= \qmod{{\cal A}(\det\Sigma^+)
\times(A^0(\Sigma^+)\setminus\{0\})}{{\cal G}}$.

The space ${\cal B}(c)^*$ has the weak homotopy type of
$K(\Z,2)\times K(H^1(X,\Z),1)$ and there is a natural isomorphism
$\nu: \Z[u]\otimes
\Lambda^*(\qmod{H_1(X,\Z)}{\rm Tors})\map H^*({\cal B}(c)^*,\Z)$.

Suppose  that $(g,b)$ is a c-good pair and fix
$\cg\in\pi_0(\qmod{{\cal C}}{\Aut(\hat P)})$ . The moduli space
${\cal W}^{\gamma}_{X,\beta} $ is  a compact manifold of dimension
$w_c:=\frac{1}{4}(c^2-2e(X)-3\sigma(X))$  for every  lift  $\gamma$
of $(g,\cg)$ and every regular value
$\beta$ of
$p_{\gamma,b}:{\cal W}^{\gamma}_{X,b}\map  b$. It can be oriented
by using the canonical complex orientation of the line bundle
$\det{}_\R({\rm index} (\Dr)) $
over
${\cal B}(c)^*$ together with a chosen orientation $\oo$ of the line
$\det(H^1(X,\R))\otimes\det(\H^2_{+,g}(X)^{\vee})$. Let $[{\cal
W}^{\gamma}_{X,\beta}  ]_{{\raisebox{-.5ex}{$\oo$}}} \in
H_{w_c}({\cal B}(c)^*,\Z)$ be the fundamental class associated  with
the choice of $\oo$.

 The
\underbar{Seiberg}-\underbar{Witten} \underbar{form}
$SW_{X,\oo}^{(g,b)}(\cg)\in \Lambda^* H^1(X,\Z)$ associated with
$(\oo,(g,b),\cg)$ is   defined by
$SW_{X,\oo}^{(g,b)}(\cg)(l_1\wedge\dots\wedge
l_r)=\left\langle\nu(l_1)\cup\dots\cup\nu(l_r)\cup
u^{\frac{w_c-r}{2}}, [{\cal W}^{\gamma}_{X,\beta}
]_{{\raisebox{-.5ex}{$\oo$}}}
\right\rangle$
for decomposable elements $l_1\wedge\dots\wedge l_r$
with  $r\equiv w_c$ (mod 2).
Here $\gamma$ lifts the pair $(g,\cg)$ and $\beta\in b$ is a  regular
value of
$p_{\gamma,b}$.
The form $SW_{X,\oo}^{(g,b)}(\cg)$ is well-defined, since the
 cohomology classes $u$, $\nu(l_i)$, as well as the trivialization of
the orientation line bundle extend  to     $\qmod{{\cal A}(\det
\Sigma^+)\times (A^0(\Sigma^+)\setminus\{0\})\times b}{{\cal G}}$,   and since
the homology class defined by $[{\cal W}^{\gamma}_{X,\beta}
]_{{\raisebox{-.5ex}{$\oo$}}}$ in this quotient  depends only on
$(g_\gamma,\cg_\gamma)$ and $b$.
Now there are two cases:

If $b_+>1$, then    $SW_{X,\oo}^{(g,b)}(\cg)$
does not depend on $(g,b)$, since the
cohomology classes $u$, $\nu(l_i)$ and the trivialization of the
orientation line bundle extend to
$\Aut(\hat P)$-invariant objects on the  quotient $\qmod{{\cal
A}(\det\Sigma^+)\times (A^0(\Sigma^+)\setminus\{0\})\times
Z^2_{\rm DR}(X)\times{\cal C}}{{\cal G}}$. Thus we may simply write
$SW_{X,\oo}(\cg)\in
\Lambda^* H^1(X,\Z)$. If $b_1=0$, then we obtain numbers
$n_{\cg}^{\ooo}$ which can be considered as refinements of the
numbers
$n_c^{\ooo}$  defined in [W]. Indeed, $n_c^{\ooo}=\sum_{\cg}
n^{\ooo}_{\cg}$, the summation being over all
$\cg\in\pi_0(\qmod{{\cal C}}{\Aut(\hat P)})$.

Suppose now that $b_+=1$. There is a natural map ${\cal M}et_X\map
\P(H^2_{\rm DR}(X))$ which sends a metric $g$ to the line
$\R[\omega_+]\subset H^2_{\rm DR}(X)$, where $\omega_+$ is any
non-trivial $g$-selfdual harmonic form. The  hyperbolic space
${\bf H}:=\{\omega\in H^2_{\rm DR}(X)|\ \omega^2=1\}
$
 has two connected components, and the choice of one of them
orients the lines $\H^2_{+,g}(X)$ for all metrics $g$. Having fixed a
component
${\bf H}_0$ of ${\bf H}$, every metric defines a unique
$g$-self-dual form
$\omega_g$ with  $[\omega_g]\in{\bf H}_0$.
\vskip 5 pt
{\sc Definition.} - {\it Let $X$ be a manifold with $b_+=1$,   and let
$c\in H^2(X,\Z)$ be characteristic. The \underbar{wall} associated
with $c$ is the hypersurface
$c^{\bot}:=\{(\omega,b)\in{\bf H}\times H^2_{\rm DR}(X)|\
(c-b)\cdot\omega=0\}$.  The connected components of ${\bf
H}\setminus c^{\bot}$ are called
\underbar{chambers} of type $c$.}
\vskip 3 pt

Notice that the walls are non-linear! Every characteristic  element
$c$ defines precisely four chambers of type $c$, namely
$C_{{\bf H}_0,\pm}:=\{(\omega,b)\in{\bf H}_0\times  H^2_{\rm
DR}(X)|\
\pm(c-b)\cdot\omega<0\}\ ,
$
 where ${\bf H}_0$ is one of the components of ${\bf H}$.  Each of
these four chambers contains pairs of the form $([\omega_g],b)$. Let
$\oo_1$ be an orientation of $H^1(X,\R)$.
\vskip 5pt
{\sc Definition.} - {\it The Seiberg-Witten invariant associated
with
$(\oo_1,{\bf H}_0,\cg)$ is the function
 $ SW_{X,(\oo_1,{\bf H}_0)}(\cg):\{\pm\}\map  \Lambda^* H^1(X,\Z)$
given by  $SW_{X,(\oo_1,{\bf H}_0)}(\cg)(\pm):=
SW^{(g,b)}_{X,\oo}(\cg)$, where
$\oo$ is the orientation defined by $(\oo_1,{\bf H}_0)$, and  $(g,b)$
is a pair   such that $([\omega_g],b)\in C_{{\bf H}_0,\pm}$.
}
\vskip 3pt
Note that, changing the orientation $\oo_1$  changes  the  invariant
by a factor
$-1$, and that $SW_{X,(\oo_1,-{\bf H}_0)}(\cg)(\pm)=
-SW_{X,(\oo_1, {\bf H}_0)}(\cg)(\mp)$.
\vskip 3pt
{\sc Remark.} -  A different approach - adapting ideas from intersection
theory to
construct "Seiberg-Witten multiplicities" - has been proposed by R. Brussee.
\vskip 3pt

Define   $u_c\in \Lambda^2\left(\qmod{H_1(X,\Z)}{\rm
Tors}\right)$ by the formula
$u_c(a,b):=\frac{1}{2}\langle c\cup a\cup b,[X]\rangle $,  for
elements $a,b\in H^1(X,\Z)$.
The following wall crossing formula generalizes results of
[W],~[KM],~[LL]\footnote{After
submitting the first version of this note, we were informed that an
expanded proof of the result in [LL] appears in an unpublished manuscript
of D. Salamon.}.
\vskip 3pt
{\sc Theorem.} -  {\it Let $b^+(X)=1$,  let $l_{\ooo_1}$
be the generator of $\Lambda^{b_1}(H^1(X,\Z))$ defined by the
orientation $\oo_1$, and let $r\geq 0$, $r\equiv w_c$ (mod 2).
For every $\lambda\in \Lambda^r\left(\qmod{H_1(X,\Z)}{\rm
Tors}\right)$  we have\\
\centerline{$ SW_{X,(\oo_1,{\bf H}_0)}(\cg)(+)(\lambda)- SW_{X,(\oo_1,{\bf
H}_0)}(\cg)(-)(\lambda)=
 \frac{ (-1)^{\left[\frac{b_1-r}{2}\right]} }{ \left[\frac{b_1-r}{2}\right]!}
 \langle  \lambda\wedge u_c
^{\left[\frac{b_1-r}{2}\right]},l_{\ooo_1} \rangle$}
 if $ r\leq \min(b_1,w_c) $, and   the difference is  0  otherwise.
 }

 \vskip 5pt
{\sc Remark.} - {\it For manifolds admitting a metric of positive
scalar curvature, the invariants are determined by Witten's vanishing result [W]
and the wall crossing formula .}
\vskip 3pt

3. {\sc Seiberg-Witten invariants of K\"ahler surfaces.} - Let $(X,g)$
be a K\"ahler surface with K\"ahler form
$\omega_g$, and let $\cg_0$ be the class of the canonical
$\Spin^c(4)$-structure of determinant
$K_X^{\vee}$ on $(X,g)$. The corresponding spinor bundles are
$\Sigma^{+}=\Lambda^{00}\oplus\Lambda^{02}$, $\Sigma^{-}=
\Lambda^{01}$ [OT1]. There is a natural bijection between classes of
$\Spin^c(4)$-structures $\cg$ of Chern class $c$ and isomorphism
classes of line bundles
$M$ with $2c_1(M)-c_1(K_X)=c$. We denote by $\cg_M$ the class
defined by a line bundle  $M$. The spinor bundles of  $\cg_M$  are the
tensor products
$\Sigma^{\pm}\otimes M$, and the map $\gamma_M:\Lambda^1_X
\map \R SU(\Sigma^+\otimes M,\Sigma^-\otimes M)$ given by
$\gamma_M(\cdot)=\gamma_0(\cdot)\otimes\id_M$ is a Clifford
map representing
$\cg_M$.
 Let  ${\cal D}ou(m)$ be the Douady space of all effective
divisors $D$ on $X$ with $c_1({\cal O}_X(D))=m$.
\vskip 5pt
{\sc Theorem.} - {\it
Let $(X,g)$ be a  connected K\"ahler surface,   and let $\cg_M$ be the
class of the $\Spin^c(4)$-structure associated to a  line bundle
$M$ with $c_1(M)=m$. Let $\beta\in A^{1,1}_{\R}$ be a   closed  form
representing the
  class $b$ such that
$ (2m-c_1(K_X)-b)\cup[\omega_g]<0$ ($>0$).
\hfill{\break}
i) \ If $c\ \not\in\  NS(X)$, then ${\cal W}_{X,\beta}^{\gamma_M}
=\emptyset$.    If $c\in NS(X)$, then there is a natural real analytic
isomorphism
${\cal W}_{X,\beta}^{\gamma_M}\simeq  {\cal
D}ou(m)$ $({\cal D}ou(c_1(K_X)-m))$.
 \hfill{\break}
ii)  ${\cal W}_{X,\beta}^{\gamma_M} $ is  smooth at a point
corresponding to
$D\in{\cal D}ou(m)$ iff
$h^0({\cal O}_D(D))=\dim_D{\cal D}ou(m)\ .$
This condition is always satisfied when $h^1({\cal O}_X)=0$.
\hfill{\break}
iii) If  ${\cal W}_{X,\beta}^{\gamma_M}$ is  smooth at a point
corresponding to $D$, then it has the expected dimension in this
point iff $h^1({\cal O}_D(D))=0$. }
\vskip 3pt
It easy to see that every  connected complex surface with $p_g>0$
is oriented diffeomorphic to a surface which possesses a
0-connected canonical divisor. This yields an easy proof of
\vskip 4pt
{\sc Corollary.} - {\it ([W]) All non-trivial Seiberg-Witten invariants
of K\"ahler surfaces with $p_g>0$ have index 0.
}
\vskip 5pt
On the other hand, for $p_g=0$,  we have
\vskip 5pt
{\sc Corollary.} - {\it Let $X$ be a surface with $p_g=0$ and $q=0$.
For all data
$({\bf H}_0,\cg)$ with $w_c\geq 0$,  we have
$SW_{X,{\bf H}_0}(\cg)(\{\pm\})=\{0,1\}$ or  $SW_{X,{\bf
H}_0}(\cg)(\{\pm\})=\{0,-1\}$. }
\vskip 3pt

{\sc Remark.} - {\it In this situation, the invariants are already
determined by their reduction modulo 2.  There exist examples of
surfaces with $p_g=0$ and
$q=0$, having infinitely many non-trivial invariants of any given
non-negative index.}
\vskip 3pt
{\sc Example.} -  Let $X=\P^2$,   let $h\in H^2(\P^2,\Z)$ be the  class
of the ample generator, and let ${\bf H}_0$ be the component of
${\bf H}=\{\pm h\}$ which contains $h$.  The classes of
$\Spin^c(4)$-structures are labelled by odd integers $c$ and the
corresponding index is $w_c=\frac{1}{4}(c^2-9)$. We find $SW_{\P^2,{\bf
H}_0}(\pm 1)\equiv 0$,
$SW_{\P^2,{\bf H}_0}(c)(+)=  1$ if $c\geq 3$, and $SW_{\P^2,{\bf
H}_0}(c)(-)= - 1$ if $c\leq 3$.
\parindent0cm
\vspace{0.3cm}\\
{{\sc References}}\vskip 8pt
{\small
[DK] Donaldson, S.; Kronheimer, P.: {\it The Geometry of
four-manifolds}, Oxford Sc. Publ. 1990.

[KM] Kronheimer, P.; Mrowka, T.: {\it The genus of embedded
surfaces in  the projective plane}, Math. Res. Lett. 1,  1994,
797-808.

[LL] Li, T.; Liu, A.: {\it General wall crossing formula}, Math. Res. Lett.
2,  1995,
797-810.

[OT1]   Okonek, Ch.; Teleman A.: {\it The Coupled Seiberg-Witten
Equations, Vortices, and Moduli Spaces of Stable Pairs},  Intern. J.
Math. Vol. 6, No. 6,  1995,
 893-910.

[OT2] Okonek, Ch.; Teleman A.: {\it Quaternionic monopoles}, Comm.
Math. Phys. (to appear).

[T] Taubes,  C.: {\it The Seiberg-Witten and Gromov invariants} Math.
Res. Lett. 2,  1995,   221-238.

[W] Witten, E.: {\it Monopoles and four-manifolds}, Math.
Res. Lett. 1,  1994,   769-796.
\\
%
\vskip-1.5pt
\vbox{
Ch. Okonek: Mathematisches Institut, Universit\"at Z\"urich,
Winterthurerstr. 190,
CH-8057 Z\"urich \\
A. Teleman: \ Mathematisches Institut, Universit\"at Z\"urich,
Winterthurerstr. 190,
CH-8057 Z\"urich \\
\hspace*{1.8cm} and Faculty of Mathematics, University of
Bucharest\\
\hspace*{1.8cm} e-mail: okonek@math.unizh.ch ;
teleman@math.unizh.ch}}

\end{document}